\documentclass[12pt, twoside, leqno]{article}
\usepackage{amsmath,amsthm}
\usepackage{amssymb,latexsym}
\usepackage{enumerate}

\newtheorem{thm}{Theorem}[section]
\newtheorem{prop}[thm]{Proposition}
\newtheorem{cor}[thm]{Corollary}
\newtheorem{lem}[thm]{Lemma}

\newtheorem{conj}[thm]{Conjecture}

\theoremstyle{definition}

\numberwithin{equation}{section}

\usepackage{mathrsfs}
\usepackage{CJK}
\usepackage{amsmath, amscd}
\usepackage{fancyhdr}
\usepackage{enumerate}

\begin{document}

\baselineskip=17pt

\title{The images of some simple derivations}
\author{ Dan Yan \footnote{ The author is supported by Scientific Research Fund of Hunan Provincial Education Department (Grant No. 21A0056), the NSF of China (Grant No. 11871241; 11601146) and the Construct Program of the Key Discipline in Hunan Province.}\\
MOE-LCSM,\\ School of Mathematics and Statistics,\\
 Hunan Normal University, Changsha 410081, China \\
\emph{E-mail:} yan-dan-hi@163.com \\
}
\date{}

\maketitle

\renewcommand{\thefootnote}{}

\renewcommand{\thefootnote}{\arabic{footnote}}
\setcounter{footnote}{0}

\begin{abstract} In the paper, we study the relation between the images of polynomial derivations and their simplicity. We prove that the images of simple Shamsuddin derivations are not Mathieu-Zhao spaces. In addition, we also show that the images of some simple derivations in dimension three are not Mathieu-Zhao spaces.
Thus, we conjecture that the images of simple derivations in dimension greater than one are not Mathieu-Zhao spaces. We also prove that locally nilpotent derivations are not simple in dimension greater than one.
\end{abstract}
 \vskip 2.5mm

\section{Introduction}

Throughout this paper, we will write $\mathbb{N}$ for the non-negative integers,
$K$ for any field with
characteristic zero and $R:=K[x,y_1,\ldots,y_n]$ for the
polynomial algebra over $K$ in $n+1$ indeterminates $x,y_1,\ldots,y_n$.
$\partial_x,~\partial_i$ will denote the derivations $\frac{\partial}{\partial x}$, $\frac{\partial}{\partial y_i}$ of $R$ for all $1\leq i\leq n$, respectively. More generally, if $s,r_1,\ldots, r_s\geq 1$ are integers and $\{x\}\bigcup\{y_{i,j}:i=1,\ldots,s,~j=1,\ldots,r_i\}$ are indeterminates over $K$, $\partial_{i,j}$ will denote the derivation $\frac{\partial}{\partial y_{i,j}}$ of $K[x,\bigcup_{i=1}^s\{y_{i,1},\ldots,y_{i,r_i}\}]$. We abbreviate $\frac{\partial g_t}{\partial y_j}$ as $g_{ty_j}$. For element $f$ of $K[x]$, we shall often use $f'$ instead of $f_x$.

A $K$-derivation $D:R\rightarrow R$ of $R$ is a $
K$-linear map such that
$$D(ab)=D(a)b+aD(b)$$
for any $a,b\in R$ and $D(c)=0$ for any $c\in K$. The set of all $K$-derivations of $R$
is denoted by $\operatorname{Der}_ K(R)$. An ideal $I$ of $R$ is called $D$-$stable$ if $D(I)\subset I$. $R$ is called
$D$-$simple$ if it has no proper nonzero $D$-$stable$ ideal. The $K$-derivation $D$ is called $simple$ if $R$ has no $D$-$stable$
ideals other than $0$ and $R$. For some examples of simple
derivations, see \cite{4}, \cite{9}, \cite{5}, \cite{2}.

A derivation $D$ of $R$ is said to be a Shamsuddin derivation if
$D=\partial_x+\sum_{i=1}^n(a_iy_i+b_i)\partial_i$ with $a_i,
b_i\in  K[x]$ for all $1\leq i\leq n$. Observe that if $D$ is such a Shamsuddin derivation of $R$, then grouping the terms that have the same $a_i$ and renaming the indeterminates $y_i$ and the polynomials $a_i, b_i$ if necessary, we can write $D$ in the following form:
$$D=\partial_x+\sum_{i=1}^s\sum_{j=1}^{r_i}(a_iy_{i,j}+b_{i,j})\partial_{i,j}$$
with $a_i, b_{i,j}\in K[x]$ for every $i$ and every $(i,j)$, $a_i\neq a_l$ for $i\neq l$. A derivation $D$ of $R$ is said to be locally nilpotent if, for each $a\in R$, there exists $m\in \mathbb{N}^*$ such that $D^m(a)=0$.

The Mathieu-Zhao space was introduced by Zhao in \cite{12} and \cite{14}, which is a natural generalization of ideals. We give the definition here for the polynomial rings. A $K$-subspace $M$ of $R$ is said to be a Mathieu-Zhao space if for any $a, b\in R$ with $a^m\in M$ for all $m\geq 1$, we have $ba^m\in M$ when $m\gg 0$.

In our paper, we prove that the images of simple Shamsuddin derivations are not Mathieu-Zhao spaces in section 2. In addition, we give a necessary and sufficient condition for $D$ to be a Mathieu-Zhao spaces, where $D=\partial_x+\sum_{i=1}^s(a(x)y_i+b_i(x))\partial_i+\sum_{j=s+1}^nb_j(x)\partial_j$. In section 3, we prove that the images of some simple derivations in dimension three are not Mathieu-Zhao spaces and locally nilpotent derivations are not simple in dimension greater than one. According to our conclusions, we make the following conjecture:

\begin{conj}
If $D$ is a simple derivation
of $K[x_1,\ldots,x_n]$ and $n\geq 2$, then $\operatorname{Im}D$ is not a Mathieu-Zhao space.
\end{conj}

\section{The images of simple Shamsuddin derivations}

\begin{lem} \label{lem2.1}
Let $D=\partial_x+\sum_{i=1}^s\sum_{j=1}^{r_i}(a_iy_{i,j}+b_{i,j})\partial_{i,j}$ be a derivation of $K[x,y_{1,1},\allowbreak\ldots,y_{s,r_s}]$ with $a_i, b_{i,j}\in K[x]$. If $D$ is simple, then $\deg a_i\geq 1$ for all $1\leq i\leq s$.
\end{lem}
\begin{proof}
It follows from Theorem 3.1 in \cite{6} that $D_i=\partial_x+\sum_{j=1}^{r_i}(a_iy_{i,j}+b_{i,j})\partial_{i,j}$ is simple for all $1\leq i\leq s$. If $a_i\in K$, then it's easy to see that the equation $z'=a_iz+\sum_{j=1}^{r_i}k_jb_{i,j}$ has solution in $K[x]$ for any $(k_1,\ldots,k_{r_i})\in K^{r_i}$. It follows from Theorem 3.2 in \cite{6} that $D_i$ is not simple, a contradiction. Hence $\deg a_i\geq 1$ for all $1\leq i\leq s$.
\end{proof}

\begin{thm}\label{thm2.2}
Let $D=\partial_x+\sum_{i=1}^s\sum_{j=1}^{r_i}(a_iy_{i,j}+b_{i,j})\partial_{i,j}$ be a derivation of $K[x,y_{1,1},\ldots,y_{s,r_s}]$ with $a_i, b_{i,j}\in K[x]$. If $D$ is simple, then $\operatorname{Im}D$ is not a Mathieu-Zhao space of $K[x,y_{1,1},\ldots,y_{s,r_s}]$.
\end{thm}
\begin{proof}
Note that $1\in \operatorname{Im}D$. If $\operatorname{Im}D$ is a Mathieu-Zhao space of $K[x,y_{1,1},\ldots,y_{s,r_s}]$, then we have $\operatorname{Im}D=K[x,y_{1,1},\ldots,y_{s,r_s}]$. We claim that $y_{1,1}\notin \operatorname{Im}D$. Suppose that $y_{1,1}\notin \operatorname{Im}D$. Then there exists $f\in K[x,y_{1,1},\ldots,y_{s,r_s}]$ such that
\begin{eqnarray}\label{eq2.1}
y_{1,1}=D(f).
\end{eqnarray}
Let $f=f^{(d)}+f^{(d-1)}+\cdots+f^{(1)}+f^{(0)}$, where $f^{(j)}$ is a polynomial of degree $j$ with respect to $y_{1,1},\ldots,y_{s,r_s}$ for $0\leq j\leq d$. Suppose that
$$f^{(d)}=\sum_{l_1+\cdots+l_s=d}c_{l_1,\ldots,l_s}y_1^{l_1}\cdots y_s^{l_s}$$
with $c_{l_1,\ldots,l_s}\in K[x]$, $y_i=(y_{i,1},\ldots,y_{i,r_i})$, $l_i=l_{i,1}+\cdots+l_{i,r_i}$ for $1\leq i\leq s$.
If the equation $a_1(x)\gamma_1+\cdots+a_s(x)\gamma_s=0$ has no positive integral solution, then the conclusion follows from Lemma \ref{lem2.1} and Theorem 3.2 in \cite{3}. Otherwise, we can assume that $\sum_{i=1}^s|m_i^{(u)}|a_i(x)=0$ with $(m_1^{(u)},\ldots,m_s^{(u)})=1$ and $m_i^{(u)}=(m_{i,1}^{(u)},\ldots,m_{i,r_i}^{(u)})$ for all $1\leq i\leq s$, $1\leq u\leq t$ for some $t\in\{1,2,\ldots,s\}$, $|m_i^{(u)}|=m_{i,1}^{(u)}+\cdots+m_{i,r_i}^{(u)}$ and the vectors $m^{(1)}:=(m_1^{(1)},\ldots,m_s^{(1)}),\ldots,m^{(t)}:=(m_1^{(t)},\ldots,m_s^{(t)})$ are linearly independent. We view the polynomials as in $K[x,y_{1,1},\allowbreak\ldots,y_{s,r_s}]$ with coefficients in $K[x]$ when we comparing the coefficients of monomials on $y_{1,1},\ldots,y_{s,r_s}$.

$(1)$ If $d\geq 2$, then we have
\begin{eqnarray}\label{eq2.2}
c_{l_1,\ldots,l_s}'+(l_1a_1(x)+\cdots+l_sa_s(x))c_{l_1,\ldots,l_s}=0
\end{eqnarray}
by comparing the coefficients of $y_1^{l_1}\cdots y_s^{l_s}$ with $l_1+\cdots+l_s=d$ of equation \eqref{eq2.1}. If $l_1a_1(x)+\cdots+l_sa_s(x)\neq 0$, then it follows from equation \eqref{eq2.2} that $c_{l_1,\ldots,l_s}=0$. If $l_1a_1(x)+\cdots+l_sa_s(x)= 0$, then $c_{l_1,\ldots,l_s}\in K$. Thus, we have
$$f^{(d)}=c_{k_1}(y_1^{m_1^{(1)}}\cdots y_s^{m_s^{(1)}})^{k_1}+\cdots+c_{k_t}(y_1^{m_1^{(t)}}\cdots y_s^{m_s^{(t)}})^{k_t}$$
with $c_{k_1}\neq 0$, where $\sum_{i=1}^sk_u|m_i^{(u)}|=d$ for all $1\leq u\leq t$. Suppose that $f^{(d-1)}=\sum_{l_1+\cdots+l_s=d-1}\tilde{c}_{l_1,\ldots,l_s}y_1^{l_1}\cdots y_s^{l_s}$ with $\tilde{c}_{l_1,\ldots,l_s}\in K[x]$, $y_i=(y_{i,1},\ldots,y_{i,r_i})$, $l_i=l_{i,1}+\cdots+l_{i,r_i}$ for $1\leq i\leq s$. Note that $(k_um_{1,1}^{(u)},\ldots,k_um_{i,j}^{(u)}-1,\ldots,k_um_{s,r_s}^{(u)})\neq (k_vm_{1,1}^{(v)},\ldots,k_vm_{\tilde{i},\tilde{j}}^{(v)}-1,\ldots,k_vm_{s,r_s}^{(v)})$ for $i\neq \tilde{i}$, $1\leq u, v\leq t$. Otherwise, $a_i(x)=a_{\tilde{i}}(x)$, a contradiction.

If $d\geq 3$, then comparing the coefficients of $y_{1,1}^{k_1m_{1,1}^{(1)}-1}\cdots y_{1,r_1}^{k_1m_{1,r_1}^{(1)}}\cdots y_{s,1}^{k_1m_{s,1}^{(1)}}\cdots \allowbreak y_{s,r_s}^{k_1m_{s,r_s}^{(1)}}$, $y_{1,1}^{k_1m_{1,1}^{(1)}} y_{1,2}^{k_1m_{1,2}^{(1)}-1}\cdots y_{s,r_s}^{k_1m_{s,r_s}^{(1)}}$, $\ldots$ , $y_{1,1}^{k_1m_{1,1}^{(1)}} y_{1,2}^{k_1m_{1,2}^{(1)}}\cdots y_{s,r_s}^{k_1m_{s,r_s}^{(1)}-1}$ of equation \eqref{eq2.1}, at least one of $\tilde{c}_{l_1,\ldots,l_s}$ satisfies the equation $z'=a_{i_0}z+\sum_{j=1}^{r_{i_0}}\hat{k}_jb_{i_0,j}$ for some $i_0\in \{1,2,\ldots,s\}$, $(\hat{k}_1,\ldots,\hat{k}_{r_{i_0}})\in K^{r_{i_0}}/\{(0,\ldots,0)\}$; $l_1+\cdots+l_s=d-1$. Since $D$ is simple, it follows from Theorem 3.1 in \cite{6} that $D_i=\partial_x+\sum_{j=1}^{r_i}(a_iy_{i,j}+b_{i,j})\partial_{i,j}$ is simple for all $1\leq i\leq s$. It follows from Theorem 3.2 in \cite{6} that the equation  $z'=a_{i_0}z+\sum_{j=1}^{r_{i_0}}\tilde{k}_jb_{i_0,j}$ does not have any solution in $K[x]$ for any $(\tilde{k}_1,\ldots,\tilde{k}_{r_{i_0}})\in K^{r_{i_0}}/\{(0,\ldots,0)\}$, a contradiction.

If $d=2$, then $f^{(2)}=c_{k_1}y_{i_1,j_1}y_{\tilde{i}_1,\tilde{j}_1}+\cdots+c_{k_t}y_{i_t,j_t}y_{\tilde{i}_t,\tilde{j}_t}$ with $c_{k_1}\neq 0$, $i_u\neq \tilde{i}_u$, $1\leq u\leq t$. Let $f^{(1)}=\sum_{i=1}^s\sum_{j=1}^{r_i}C_{i,j}y_{i,j}$ with $C_{i,j}\in K[x]$ for $1\leq j\leq r_i$, $1\leq i\leq s$. If $y_{1,1}\notin \{y_{i_1,j_1},\ldots,y_{i_t,j_t}, y_{\tilde{i}_1,\tilde{j}_1},\ldots,y_{\tilde{i}_t,\tilde{j}_t}\}$, then we have
\begin{eqnarray}\label{eq2.3}
C_{1,1}'+a_1(x)C_{1,1}=1
\end{eqnarray}
by comparing the coefficients of $y_{1,1}$ of equation \eqref{eq2.1}. It follows from Lemma \ref{lem2.1} that $\deg a_1(x)\geq 1$. Then we have $C_{1,1}=0$ by comparing the degree of $x$ of equation \eqref{eq2.3}, which contradicts equation \eqref{eq2.3}. If $y_{1,1}\in \{y_{i_1,j_1},\ldots,y_{i_t,j_t}, y_{\tilde{i}_1,\tilde{j}_1},\ldots,y_{\tilde{i}_t,\tilde{j}_t}\}$, then we can assume that
$$f^{(2)}=c_{k_1}y_{1,1}y_{i_1,j_1}+\cdots+c_{k_t}y_{i_t,j_t}y_{\tilde{i}_t,\tilde{j}_t}.$$
Hence
\begin{eqnarray}\label{eq2.4}
C_{i_1,j_1}'-a_1(x)C_{i_1,j_1}+c_{k_1}b_{1,1}+c_{k_2}b_{1,2}+\cdots+c_{k_{\tilde{t}_1}}b_{1,\tilde{t}_1}=0
\end{eqnarray}
by comparing the coefficients of $y_{i_1,j_1}$ of equation \eqref{eq2.1}, where $1\leq \tilde{t}_1\leq \operatorname{min}\{t,r_1\}$. Since $D$ is simple, it follows from the arguments of $d=3$ that $C_{i_1,j_1}'-a_1(x)C_{i_1,j_1}+c_{k_1}b_{1,1}+c_{k_2}b_{1,2}+\cdots+c_{k_{\tilde{t}_1}}b_{1,\tilde{t}_1}\neq 0$, a contradiction.

$(2)$ If $d=1$, then we can assume that $f^{(1)}=\sum_{i=1}^s\sum_{j=1}^{r_i}C_{i,j}y_{i,j}$ with $C_{i,j}\in K[x]$ for $1\leq j\leq r_i$, $1\leq i\leq s$. Thus, we have
\begin{eqnarray}\label{eq2.5}
C_{1,1}'+a_1(x)C_{1,1}=1
\end{eqnarray}
by comparing the coefficients of $y_{1,1}$ of equation \eqref{eq2.1}. Since $\deg a_1(x)\geq 1$, we have $C_{1,1}=0$ by comparing the degree of $x$ of equation \eqref{eq2.5}, which contradicts equation \eqref{eq2.5}.

$(3)$ If $d=0$, then $D(f)\in K[x]$. Clearly, $y_{1,1}\notin \operatorname{Im}D$. Thus, the conclusion follows.
\end{proof}

\begin{cor}\label{cor2.3}
Let $D=\partial_x+\sum_{i=1}^s\sum_{j=1}^{r_i}(a_i(x)y_{i,j}+b_{i,j})\partial_{i,j}$ be a derivation of $K[x,y_{1,1},\ldots,y_{s,r_s}]$ with $a_i, b_{i,j}\in K[x]$. If $\deg a_i>\deg b_{i,j}$ and $b_{i,1},\ldots,b_{i,r_i}$ are linearly independent over $K$ for every $i\in \{1,2,\ldots,s\}$ and every $j\in \{1,2,\ldots,r_i\}$, then $ImD$ is not a Mathieu-Zhao space of $K[x,y_{1,1},\ldots,y_{s,r_s}]$.
\end{cor}
\begin{proof}
It follows from Corollary 3.3 in \cite{6} that $D$ is simple. Then the conclusion follows from Theorem \ref{thm2.2}.
\end{proof}

\begin{prop}\label{prop2.4}
Let $D=\partial_x+\sum_{i=1}^k(a(x)y_i+b_i(x))\partial_i+\sum_{j=k+1}^nb_j(x)\partial_j$ be a derivation of $K[x,y_1,\ldots,y_n]$ with $a(x), b_i(x)\in K[x]$ for all $1\leq i\leq n$, $k\geq 1$. Then $\operatorname{Im}D$ is a Mathieu-Zhao space of $K[x,y_1,\ldots,y_n]$ iff $a(x)\in K$.
\end{prop}
\begin{proof}
$``\Leftarrow"$ If $a(x)\in K$, then it follows from Example 9.3.2 in \cite{10} that $D$ is locally finite. Since $1\in \operatorname{Im}D$, it follows from Proposition 1.4 in \cite{8} that $\operatorname{Im}D$ is a Mathieu-Zhao space of $K[x,y_1,\ldots,y_n]$.

$``\Rightarrow"$ Since $1\in \operatorname{Im}D$ and $\operatorname{Im}D$ is a Mathieu-Zhao space of $K[x,y_1,\ldots,y_n]$, we have $\operatorname{Im}D=K[x,y_1,\ldots,y_n]$. Suppose that $\deg a(x)\geq 1$. We claim $y_1\notin \operatorname{Im}D$. If $y_1\in \operatorname{Im}D$, then there exists $f\in K[x,y_1,\ldots,y_n]$ such that
\begin{eqnarray}\label{eq2.6}
y_1=D(f).
\end{eqnarray}
Let $f=f^{(d)}+\cdots+f^{(1)}+f^{(0)}$, where $f^{(\tilde{k})}$ is a polynomial of degree $k$ with respect to $y_1,\ldots,y_n$ for $0\leq \tilde{k}\leq d$. Suppose that $f^{(d)}=\sum_{l_1+\cdots+l_n=d}c_{l_1,\ldots,l_n}^{(0)}y_1^{l_1}\cdots y_n^{l_n}$ with $c_{l_1,\ldots,l_n}^{(0)}\in K[x]$.

If $d\geq 2$, then we have
\begin{eqnarray}\label{eq2.7}
c_{l_1,\ldots,l_n}^{(0)'}+(l_1+\cdots+l_k)a(x)c_{l_1,\ldots,l_n}^{(0)}=0
\end{eqnarray}
by comparing the coefficients of $y_1^{l_1}\cdots y_n^{l_n}$ with $l_1+\cdots+l_n=d$ of equation \eqref{eq2.6}.
If $l_1+\cdots+l_k\neq 0$, then $c_{l_1,\ldots,l_n}^{(0)}=0$ for all $l_1+\cdots+l_n=d$ by comparing the degree of $x$ of equation \eqref{eq2.7}. Thus,
$$f^{(d)}=\sum_{l_{k+1}+\cdots+l_n=d}c_{l_{k+1},\ldots,l_n}^{(0)}y_{k+1}^{l_{k+1}}\cdots y_n^{l_n}.$$
Let $f^{(d-1)}=\sum_{l_1+\cdots+l_n=d-1}c_{l_1,\ldots,l_n}^{(1)}y_1^{l_1}\cdots y_n^{l_n}$ with $c_{l_1,\ldots,l_n}^{(1)}\in K[x]$. If $d-1\geq 2$, then we have
\begin{eqnarray}\label{eq2.8}
\nonumber
\sum_{l_1+\cdots+l_n=d-1}[c_{l_1,\ldots,l_n}^{(1)'}+(l_1+\cdots+l_k)a(x)c_{l_1,\ldots,l_n}^{(1)}]y_1^{l_1}\cdots y_n^{l_n}+\\
\sum_{l_{k+1}+\cdots+l_n=d}\sum_{j=k+1}^nl_jb_jc_{l_{k+1},\ldots,l_n}^{(0)}y_{k+1}^{l_{k+1}}\cdots y_j^{l_j-1}\cdots y_n^{l_n}=0
\end{eqnarray}
by considering the part of degree $d-1$ of equation \eqref{eq2.6} with respect to $y_1,\ldots,y_n$.
If $l_1+\cdots+l_k\neq 0$, then we have
\begin{eqnarray}\label{eq2.9}
c_{l_1,\ldots,l_n}^{(1)'}+(l_1+\cdots+l_k)a(x)c_{l_1,\ldots,l_n}^{(1)}=0
\end{eqnarray}
by comparing the coefficients of $y_1^{l_1}\cdots y_n^{l_n}$ with $l_1+\cdots+l_n=d-1$ of equation \eqref{eq2.8}. Thus, we have $c_{l_1,\ldots,l_n}^{(1)}=0$ for all $l_1+\cdots+l_k\neq 0$ by comparing the degree of $x$ of equation \eqref{eq2.9}. That is,
$$f^{(d-1)}=\sum_{l_{k+1}+\cdots+l_n=d-1}c_{l_{k+1},\ldots,l_n}^{(1)}y_{k+1}^{l_{k+1}}\cdots y_n^{l_n}.$$
Hence we have $f^{(t)}=\sum_{l_{k+1}+\cdots+l_n=t}c_{l_{k+1},\ldots,l_n}^{(d-t)}y_{k+1}^{l_{k+1}}\cdots y_n^{l_n}$ by considering the part of degree $t$ of equation \eqref{eq2.6} with respect to $y_1,\ldots,y_n$ for all $2\leq t\leq d$. Suppose that $f^{(1)}=c_1(x)y_1+\cdots+c_n(x)y_n$. Then we have
\begin{eqnarray}\label{eq2.10}
c_1'(x)+a(x)c_1(x)=1
\end{eqnarray}
by comparing the coefficients of $y_1$ of equation \eqref{eq2.6}. Since $\deg a(x)\geq 1$, we have $c_1(x)=0$ by comparing the degree of $x$ of equation \eqref{eq2.10}, which contradicts equation \eqref{eq2.10}.
\end{proof}

\begin{lem}\label{lem2.5}
Let $D=\partial_x+\sum_{j=1}^{r_i}(a_i(x)y_{i,j}+b_{i,j}(x))\partial_{i,j}+\sum_{l\neq i, l=1}^s\sum_{j=1}^{r_l}(a_l(x)y_{l,j})\partial_{l,j}$ be a derivation of $K[x,y_{1,1},\ldots,y_{s,r_s}]$ and $D_i=\partial_x+\sum_{j=1}^{r_i}(a_i(x)y_{i,j}+b_{i,j}(x))\partial_{i,j}$ a derivation of $K[x,y_{i,1},\ldots,y_{i,r_i}]$ with $a_i(x), a_l(x), b_{i,j}(x)\in K[x]$ for $1\leq i\leq s$. If $\operatorname{Im}D$ is a Mathieu-Zhao space of $K[x,y_{1,1},\ldots,y_{s,r_s}]$, then $\operatorname{Im}D_i$ is a Mathieu-Zhao space of $K[x,y_{i,1},\ldots,y_{i,r_i}]$ for $1\leq i\leq s$.
\end{lem}
\begin{proof}
Since $1\in \operatorname{Im}D$ and $\operatorname{Im}D$ is a Mathieu-Zhao space of $K[x,y_{1,1},\ldots,y_{s,r_s}]$, we have $\operatorname{Im}D=K[x,y_{1,1},\ldots,y_{s,r_s}]$. Thus, for any $x^{q_0}y_{i,1}^{q_1}\cdots y_{i,r_i}^{q_i}\in K[x,y_{i,1},\ldots,\allowbreak y_{i,r_i}]$, there exists $g_{q_0,\ldots,q_i}\in K[x,y_{1,1},\ldots,y_{s,r_s}]$ such that
\begin{eqnarray}\label{eq2.11}
x^{q_0}y_{i,1}^{q_1}\cdots y_{i,r_i}^{q_i}=D(g_{q_0,\ldots,q_i}).
\end{eqnarray}
Let $y_{l,j}=0$ for all $1\leq j\leq r_l$, $1\leq l\leq s$, $l\neq i$. Then equation \eqref{eq2.11} has the following form:
$$x^{q_0}y_{i,1}^{q_1}\cdots y_{i,r_i}^{q_i}=D_i(g_{q_0,\ldots,q_i}(x,0,\ldots,0,y_{i,1},\ldots,y_{i,r_i},0,\ldots,0)).$$
Hence we have $\operatorname{Im}D_i=K[x,y_{i,1},\ldots,y_{i,r_i}]$. That is, $\operatorname{Im}D_i$ is a Mathieu-Zhao space of $K[x,y_{i,1},\ldots,y_{i,r_i}]$ for $1\leq i\leq s$.
\end{proof}

\begin{cor}\label{cor2.6}
Let $D=\partial_x+\sum_{j=1}^{r_i}(a_i(x)y_{i,j}+b_{i,j}(x))\partial_{i,j}+\sum_{l\neq i, l=1}^s\sum_{j=1}^{r_l}(a_l(x)y_{l,j})\partial_{l,j}$ be a derivation of $K[x,y_{1,1},\ldots,y_{s,r_s}]$ with $a_i(x), a_l(x), b_{i,j}(x)\in K[x]$. If $\operatorname{Im}D$ is a Mathieu-Zhao space of $K[x,y_{1,1},\ldots,y_{s,r_s}]$, then $a_i(x)\in K$.
\end{cor}
\begin{proof}
The conclusion follows from Lemma \ref{lem2.5} and Corollary 3.5 in \cite{3}.
\end{proof}

\section{Images of some simple derivations}

\begin{prop}\label{prop3.1}
Let $D=\partial_{x_1}+\sum_{i=2}^n(a_ix_i+b_i)\partial_{x_i}$ be a derivation of $K[x_1,\ldots,x_n]$ with $a_i,b_i\in K[x_1,\ldots,x_{i-1}]$. Then we have the following statements:

$(1)$ If $\deg_{x_{n-1}}a_n\geq 1$, then $\operatorname{Im}D$ is not a Mathieu-Zhao space of $K[x_1,\ldots,x_n]$.

$(2)$ If $\deg_{x_{n-1}}a_n=0$ and $\deg_{x_{n-2}}(n_1a_{n-1}+n_2a_n)\geq 1$ for any $(n_1,n_2)\in \mathbb{N}^2\setminus \{(0,0)\}$, then $\operatorname{Im}D$ is not a Mathieu-Zhao space of $K[x_1,\ldots,x_n]$.
\end{prop}
\begin{proof}
Note that $1\in \operatorname{Im}D$. If $\operatorname{Im}D$ is a Mathieu-Zhao space of $K[x_1,\ldots,x_n]$, then we have $\operatorname{Im}D=K[x_1,\ldots,x_n]$. We claim $x_n\notin \operatorname{Im}D$. If $x_n\in \operatorname{Im}D$, then there exists $f\in K[x_1,\ldots,x_n]$ such that
\begin{eqnarray}\label{eq3.1}
x_n=D(f).
\end{eqnarray}
We view the polynomials as in $K[x_1,\ldots,x_i]$ with coefficients in $K[x_1,\ldots,x_{i-1}]$ when we comparing the coefficients of monomials on $x_i$ in the following arguments for $i\in \{2,\ldots,,n\}$.
Let $f=f_dx_n^d+\cdots+f_1x_n+f_0$ with $f_d\neq 0$, $f_j\in K[x_1,\ldots,x_{n-1}]$ for $0\leq j\leq d$. If $d\geq 1$, then we have
\begin{eqnarray}\label{eq3.2}
f_{dx_1}+\sum_{i=2}^{n-1}(a_ix_i+b_i)f_{dx_i}+da_nf_d=c
\end{eqnarray}
by comparing the coefficients of $x_n^d$ of equation \eqref{eq3.1}, where $c=0$ or 1.

$(1)$ If $\deg_{x_{n-1}}a_n\geq 1$, then we have $f_d=0$ by comparing the degree of $x_{n-1}$ of equation \eqref{eq3.2}, which is a contradiction. Thus, we have $d=0$. Since $D(f)\in K[x_1,\ldots,x_{n-1}]$, we have $x_n\notin \operatorname{Im}D$. Thus, $\operatorname{Im}D$ is not a Mathieu-Zhao space of $K[x_1,\ldots,x_n]$.

$(2)$ If $\deg_{x_{n-1}}a_n=0$, then we assume that $f_d=f_d^{(m)}x_{n-1}^m+\cdots+f_d^{(1)}x_{n-1}+f_d^{(0)}$ with $f_d^{(m)}\neq 0$, $f_d^{(l)}\in K[x_1,\ldots,x_{n-2}]$ for $0\leq l\leq m$, we have
\begin{eqnarray}\label{eq3.3}
f_{dx_1}^{(m)}+\sum_{i=2}^{n-2}(a_ix_i+b_i)f_{dx_i}^{(m)}+ma_{n-1}f_d^{(m)}+da_nf_d^{(m)}=c
\end{eqnarray}
by comparing the coefficients of $x_{n-1}^m$ of equation \eqref{eq3.2}, where $c=0$ or 1. Since $d\geq 1$ and $\deg_{x_{n-2}}(ma_{n-1}+da_n)\geq 1$, we have $f_d^{(m)}=0$ by comparing the degree of $x_{n-1}$ of equation \eqref{eq3.3}, which is a contradiction. Hence we have $x_n\notin \operatorname{Im}D$. Then the conclusion follows.
\end{proof}

\begin{cor}\label{cor3.2}
Let $D=\partial_{x_1}+\sum_{i=2}^n(a_ix_i+b_i)\partial_{x_i}$ be a derivation of $K[x_1,\ldots,x_n]$ with $a_i,b_i\in K[x_1,\ldots,x_{i-1}]$. If $D$ is simple, then $\deg a_2\geq 1$.
\end{cor}
\begin{proof}
Let $\tilde{D}=\partial_{x_1}+(a_2x_2+b_2)\partial_{x_2}$. Then $\tilde{D}$ is simple. Hence the conclusion follows from Lemma \ref{lem2.1}.
\end{proof}

\begin{prop}\label{prop3.3}
Let $D=\partial_{x_1}+\sum_{i=2}^3(a_ix_i+b_i)\partial_{x_i}$ be a derivation of $K[x_1,x_2,x_3]$ with $a_i,b_i\in K[x_1,\ldots,x_{i-1}]$ and $a_3=0$. If $D$ is simple, then $\operatorname{Im}D$ is not a Mathieu-Zhao space of $K[x_1,x_2,x_3]$.
\end{prop}
\begin{proof}
Note that $1\in \operatorname{Im}D$. If $\operatorname{Im}D$ is a Mathieu-Zhao space of $K[x_1,x_2,x_3]$, then $\operatorname{Im}D=K[x_1,x_2,x_3]$. Suppose that $\deg_{x_2}b_3=u$. We claim $x_2^{u+1}\notin \operatorname{Im}D$ if $D$ is simple. Suppose that $x_2^{u+1}\in \operatorname{Im}D$. Then there exists $f\in K[x_1,x_2,x_3]$ such that
\begin{eqnarray}\label{eq3.4}
x_2^{u+1}=D(f)
\end{eqnarray}
We view the polynomials as in $K[x_1,\ldots,x_i]$ with coefficients in $K[x_1,\ldots,x_{i-1}]$ for $i\in \{2,3\}$ when we comparing the coefficients of monomials on $x_i$ in the following arguments.
Let $f=f_dx_3^d+\cdots+f_1x_3+f_0$ with $f_d\neq 0$, $f_j\in K[x_1,x_2]$ for $0\leq j\leq d$. If $d\geq 1$, then we have
$$f_{dx_1}+(a_2x_2+b_2)f_{dx_2}=0$$
by comparing the coefficients of $x_3^d$ of equation \eqref{eq3.4}. That is, $D(f_d)=0$. Since $D$ is simple, we have $f_d\in K^*$. If $d\geq 2$, then we have
$$f_{(d-1)x_1}+(a_2x_2+b_2)f_{(d-1)x_2}+b_3\cdot df_d=0$$
by comparing the coefficients of $x_3^{d-1}$ of equation \eqref{eq3.4}. Then $-(df_d)^{-1}f_{d-1}$ is a solution of $\tilde{D}(z)=b_3$, where $\tilde{D}=\partial_{x_1}+(a_2x_2+b_2)\partial_{x_2}$, which contradicts that $D$ is simple.
Hence we have $d\leq 1$.

$(1)$ If $d=1$, then we have $f=f_1x_3+f_0$ with $f_1\in K^*$ and equation \eqref{eq3.4} has the following form
\begin{eqnarray}\label{eq3.5}
f_{0x_1}+(a_2x_2+b_2)f_{0x_2}+f_1\cdot b_3=x_2^{u+1}
\end{eqnarray}
Since $\deg_{x_2}b_3=u$, we have $\deg_{x_2}f_0\geq u+1$. Let $f_0=f_0^{(t)}x_2^t+\cdots+f_0^{(1)}x_2+f_0^{(0)}$ with $f_0^{(t)}\neq 0$, $f_0^{(j)}\in K[x_1]$ for $0\leq j\leq t$. If $t>u+1\geq 1$, then we have
\begin{eqnarray}\label{eq3.6}
f_0^{(t)'}+ta_2f_0^{(t)}=0
\end{eqnarray}
by comparing the coefficients of $x_2^t$ of equation \eqref{eq3.5}. It follows from Corollary \ref{cor3.2} that $\deg a_2\geq 1$. Then we have $f_0^{(t)}=0$ by comparing the degree of $x_1$ of equation \eqref{eq3.6}, which is a contradiction. Hence we have $t=u+1$. Then we have
\begin{equation}\label{eq3.7}
f_0^{(u+1)'}+ta_2f_0^{(u+1)}=1
\end{equation}
by comparing the coefficients of $x_2^{u+1}$ of equation \eqref{eq3.5}. We have a contradiction by comparing the degree of $x_1$ of equation \eqref{eq3.7}.

$(2)$ If $d=0$, then we have $f=f_0$ and equation \eqref{eq3.4} has the following form:
\begin{equation}\label{eq3.8}
f_{0x_1}+(a_2x_2+b_2)f_{0x_2}=x_2^{u+1}
\end{equation}
We have a contradiction by following the arguments of case $(1)$. Thus, we have $x_2^{u+1}\notin \operatorname{Im}D$. Hence $\operatorname{Im}D$ is not a Mathieu-Zhao space of $K[x_1,x_2,x_3]$.
\end{proof}

\begin{thm}\label{thm3.4}
Let $D=\partial_{x_1}+\sum_{i=2}^3(a_ix_i+b_i)\partial_{x_i}$ be a derivation of $K[x_1,x_2,x_3]$ with $a_i,b_i\in K[x_1,\ldots,x_{i-1}]$. If $D$ is simple, then $\operatorname{Im}D$ is not a Mathieu-Zhao space of $K[x_1,x_2,x_3]$.
\end{thm}
\begin{proof}
If $a_3=0$, then it follows from Proposition \ref{prop3.3} that $\operatorname{Im}D$ is not a Mathieu-Zhao space of $K[x_1,x_2,x_3]$. If $\deg_{x_2}a_3\geq 1$, then it follows from Proposition \ref{prop3.1} that $\operatorname{Im}D$ is not a Mathieu-Zhao space of $K[x_1,x_2,x_3]$. Hence we can assume that $a_3\neq 0$ and $\deg_{x_2}a_3=0$. Note that $1\in \operatorname{Im}D$. If $\operatorname{Im}D$ is a Mathieu-Zhao space of $K[x_1,x_2,x_3]$, then we have $\operatorname{Im}D=K[x_1,x_2,x_3]$. We claim $x_2\notin \operatorname{Im}D$ if $D$ is simple.
If $x_2\in \operatorname{Im}D$, then there exists $f\in K[x_1,x_2,x_3]$ such that
\begin{eqnarray}\label{eq3.9}
x_2=D(f).
\end{eqnarray}
We view the polynomials as in $K[x_1,\ldots,x_i]$ with coefficients in $K[x_1,\ldots,x_{i-1}]$ for $i\in \{2,3\}$ when we comparing the coefficients of monomials on $x_i$ in the following arguments.
Let $f=f_dx_3^d+\cdots+f_1x_3+f_0$ with $f_d\neq 0$, $f_j\in K[x_1,x_2]$ for $0\leq j\leq d$.

If $d\geq 1$, then we have
\begin{eqnarray}\label{eq3.10}
f_{dx_1}+(a_2x_2+b_2)f_{dx_2}+df_d\cdot a_3=0
\end{eqnarray}
by comparing the coefficients of $x_3^d$ of equation \eqref{eq3.9}. Let $f_d=f_d^{(m)}x_2^m+\cdots+f_d^{(1)}x_2+f_d^{(0)}$ with $f_d^{(m)}\neq 0$, $f_d^{(k)}\in K[x_1]$ for $0\leq k\leq d$. Then we have
\begin{eqnarray}\label{eq3.11}
f_{dx_1}^{(m)}+(ma_2+da_3)f_d^{(m)}=0
\end{eqnarray}
by comparing the coefficients of $x_2^m$ of equation \eqref{eq3.10}. Then we have $f_d^{(m)}\in K^*$ and $da_3=-ma_2$ by comparing the degree of $x_1$ of equation \eqref{eq3.11}.

If $m\geq 1$, then we have
$$f_{dx_1}^{(m-1)}+mf_d^{(m)}b_2+(m-1)f_d^{(m-1)}a_2+da_3\cdot f_d^{(m-1)}=0$$
by comparing the coefficients of $x_2^{m-1}$ of equation \eqref{eq3.10}. That is,
$$f_d^{(m-1)'}=a_2f_d^{(m-1)}-mf_d^{(m)}b_2.$$
Thus, $-(mf_d^{(m)})^{-1}f_d^{(m-1)}$ is a solution of the equation $z'=a_2z+b_2$. Let $\tilde{D}=\partial_{x_1}+(a_2x_2+b_2)\partial_{x_2}$. Then it follows from Theorem 13.2.1 in \cite{10} that $\tilde{D}$ is not simple. Hence $D$ is not simple, which is a contradiction. Whence we have $m=0$. That is, $f_d=f_d^{(0)}$. Then equation \eqref{eq3.11} has the following form:
$$f_d^{(0)'}+da_3f_d^{(0)}=0.$$
Since $a_3\neq 0$, we have $f_d^{(0)}=0$ by comparing the degree of $x_1$ of the above equation, which is a contradiction. Hence we have $d=0$. Then equation \eqref{eq3.9} has the following form:
\begin{eqnarray}\label{eq3.12}
f_{0x_1}+(a_2x_2+b_2)f_{0x_2}=x_2.
\end{eqnarray}
Let $f_0=f_0^{(t)}x_2^t+\cdots+f_0^{(1)}x_2+f_0^{(0)}$ with $f_0^{(k)}\in K[x_1]$ for $0\leq k\leq t$. Then we have $t=1$ and
\begin{eqnarray}\label{eq3.13}
f_0^{(1)'}+a_2f_0^{(1)}=1
\end{eqnarray}
by comparing the coefficients of $x_2$ of equation \eqref{eq3.12}. Since $\deg a_2\geq 1$, we have a contradiction by comparing the degree of $x_1$ of equation \eqref{eq3.13}. Hence we have $x_2\notin \operatorname{Im}D$, Whence  $\operatorname{Im}D$ is not a Mathieu-Zhao space of $K[x_1,x_2,x_3]$.
\end{proof}

\begin{prop}\label{prop3.5}
Let $D=\partial_{x_1}+(a_2x_2+b_2)\partial_{x_2}+b_3\partial_{x_3}$ be a derivation of $K[x_1,x_2,x_3]$ with $a_2,b_2\in K[x_1]$ and $b_3=b_3^{(v)}x_2^v+\cdots+b_3^{(1)}x_2+b_3^{(0)}$ with $b_3^{(v)}\neq 0$, $b_3^{(j)}\in K[x_1]$ for $0\leq j\leq v$. If $\deg a_2>\deg b_3^{(v)}$, then $\operatorname{Im}D$ is not a Mathieu-Zhao space of $K[x_1,x_2,x_3]$.
\end{prop}
\begin{proof}
Note that $1\in \operatorname{Im}D$. If $\operatorname{Im}D$ is a Mathieu-Zhao space of $K[x_1,x_2,x_3]$, then we have $\operatorname{Im}D=K[x_1,x_2,x_3]$. We claim $x_2^{v+1}\notin \operatorname{Im}D$ if $\deg a_2>\deg b_3^{(v)}$. Suppose that $x_2^{v+1}\in \operatorname{Im}D$. Then there exists $f\in K[x_1,x_2,x_3]$ such that
\begin{eqnarray}\label{eq3.14}
x_2^{v+1}=D(f).
\end{eqnarray}
We view the polynomials as in $K[x_1,\ldots,x_i]$ with coefficients in $K[x_1,\ldots,x_{i-1}]$ when we comparing the coefficients of monomials on $x_i$ for $i\in \{2,3\}$ in the following arguments. Let $f=f_dx_3^d+\cdots+f_1x_3+f_0$ with $f_d\neq 0$, $f_j\in K[x_1,x_2]$ for $0\leq j\leq d$. If $d\geq 1$, then we have the following equations:
\begin{eqnarray}\label{eq3.15}
f_{dx_1}+(a_2x_2+b_2)f_{dx_2}=0
\end{eqnarray}
by comparing the coefficients of $x_3^d$ of equation \eqref{eq3.14}.
Let $f_d=f_d^{(t)}x_2^t+\cdots+f_d^{(1)}x_2+f_d^{(0)}$ with $f_d^{(t)}\neq 0$, $f_d^{(j)}\in K[x_1]$ for $0\leq j\leq t$. If $t\geq 1$, then we have
\begin{eqnarray}\label{eq3.16}
f_d^{(t)'}+ta_2f_d^{(t)}=0
\end{eqnarray}
by comparing the coefficients of $x_2^t$ of equation \eqref{eq3.15}. Since $\deg a_2\geq 1$, we have $f_d^{(t)}=0$ by comparing the degree of $x_1$ of equation \eqref{eq3.16}, which is a contradiction. Hence we have $t=0$. It follows from equation \eqref{eq3.15} that $f_d^{(0)'}=0$. That is, $f_d=f_d^{(0)}\in K^*$.

If $d\geq 2$, then we have
\begin{eqnarray}\label{eq3.17}
f_{(d-1)x_1}+(a_2x_2+b_2)f_{(d-1)x_2}+df_d\cdot b_3=0
\end{eqnarray}
by comparing the coefficients of $x_3^{d-1}$ of equation \eqref{eq3.14}.
If $v=0$, then it follows from Proposition \ref{prop2.4} that $\operatorname{Im}D$ is not a Mathieu-Zhao space of $K[x_1,x_2,x_3]$. We can assume that $v\geq 1$. Let $f_{d-1}=f_{d-1}^{(s)}x_2^s+\cdots+f_{d-1}^{(1)}x_2+f_{d-1}^{(0)}$ with $f_{d-1}^{(s)}\neq 0$, $f_{d-1}^{(j)}\in K[x_1]$ for $0\leq j\leq s$. It follows from equation \eqref{eq3.17} that
$s\geq v$. If $s>v$, then we have
\begin{eqnarray}\label{eq3.18}
f_{d-1}^{(s)'}+sa_2f_{d-1}^{(s)}=0
\end{eqnarray}
by comparing the coefficients of $x_3^{d-1}$ of equation \eqref{eq3.17}. Since $\deg a_2\geq 1$, we have $f_{d-1}^{(s)}=0$ by comparing the degree of $x_1$ of equation \eqref{eq3.18}, which is a contradiction. Hence we have $s=v$. Then we have
\begin{eqnarray}\label{eq3.19}
f_{d-1}^{(v)'}+sa_2f_{d-1}^{(v)}+df_d\cdot b_3^{(v)}=0
\end{eqnarray}
by comparing the coefficients of $x_2^v$ of equation \eqref{eq3.17}. Since $\deg a_2> \deg b_3^{(v)}$, we have a contradiction by comparing the degree of $x_1$ of equation \eqref{eq3.19}. Thus, we have $d\leq 1$. That is, $f=f_1x_3+f_0$ with $f_1\in K$, $f_0\in K[x_1,x_2]$. Then equation \eqref{eq3.14} has the following form:
\begin{eqnarray}\label{eq3.20}
f_{0x_1}+(a_2x_2+b_2)f_{0x_2}+f_1\cdot b_3=x_2^{v+1}.
\end{eqnarray}
We have a contradiction by following the arguments of case $(1)$ of Proposition \ref{prop3.3}. Hence we have $x_2^{v+1}\notin \operatorname{Im}D$. That is, $\operatorname{Im}D$ is not a Mathieu-Zhao space of $K[x_1,x_2,x_3]$.
\end{proof}

\begin{cor}\label{cor3.6}
Let $D=\partial_{x_1}+(a_2x_2+b_2)\partial_{x_2}+b_3\partial_{x_3}$ be a derivation of $K[x_1,x_2,x_3]$ with $a_2,b_2\in K[x_1]$ and $b_3=b_3^{(v)}x_2^v+\cdots+b_3^{(1)}x_2+b_3^{(0)}$ with $b_3^{(v)}\in K^*$, $b_3^{(j)}\in K[x_1]$ for $0\leq j\leq v$. Then $\operatorname{Im}D$ is a Mathieu-Zhao space of $K[x_1,x_2,x_3]$ iff $a_2\in K$.
\end{cor}
\begin{proof}
$``\Leftarrow"$ If $a_2\in K$, then $D$ is locally finite. Since $1\in \operatorname{Im}D$, the conclusion follows from Proposition 1.4 in \cite{8}.

$``\Rightarrow"$ If $\deg a_2\geq 1$, then $\deg a_2>\deg b_3^{(v)}$. It follows from Proposition \ref{prop3.5} that $\operatorname{Im}D$ is not a Mathieu-Zhao space of $K[x_1,x_2,x_3]$. Then the conclusion follows.
\end{proof}

\begin{prop}\label{prop3.7}
Let $D$ be a derivation of $K[x_1,x_2,\ldots,x_n]$ for $n\geq 2$. If $D$ is locally nilpotent, then $D$ is not simple.
\end{prop}
\begin{proof}
If $K[x_1,x_2,\ldots,x_n]^D\neq K$, then $K[x_1,x_2,\ldots,x_n]^D$ is a $D$-stable ideal. If $K[x_1,x_2,\ldots,x_n]^D=(1)$, then $D=0$. Clearly, $D$ is not simple. If $K[x_1,x_2,\ldots,x_n]^D\allowbreak \neq (1)$, then $K[x_1,x_2,\ldots,x_n]^D$ is a proper $D$-stable ideal. Thus, $D$ is not simple. We can assume that $K[x_1,x_2,\ldots,x_n]^D=K$. It follows from Lemma 8 in \cite{15} or Theorem 2.1 in \cite{16} that $D(f)=c\cdot J(f_1,\ldots,f_{n-1},f)$ for some $c\in K^*$, where $J(f_1,\ldots,f_{n-1},f)$ is the Jacobian of $f_1,\ldots,f_{n-1},f$. Note that $\deg f_i\geq 1$ for all $1\leq i\leq n-1$. Otherwise, $D=0$. Clearly, $D$ is not simple. Then we have $D(f_1)=0$. That is, the ideal $(f_1)$ is a proper $D$-stable ideal. Hence $D$ is not simple.
\end{proof}

\end{document}